\documentclass[12pt, leqno]{article}
\usepackage{graphicx, mathrsfs,amssymb,amsfonts,amsbsy,amsthm,latexsym,amsmath}
\numberwithin{equation}{section}
\begin{document}

\newcommand{\Dh}{\hbox{\bf D}}
\newcommand{\gh}{\hbox{\bf g}}
\newcommand{\fh}{\hbox{\bf f}}
\newcommand{\uh}{\hbox{\bf u}}
\newcommand{\vh}{\hbox{\bf v}}
\newcommand{\zh}{\hbox{\bf z}}
\newcommand{\Xh}{\hbox{\bf X}}
\newcommand{\Fh}{\hbox{\bf F}}

\title{\vskip.5cm
Differentiation matrices for meromorphic functions
\vskip.5cm}
\author{Rafael G. Campos and Claudio Meneses\\ 
Escuela de Ciencias F\'{\i}sico-Matem\'aticas\\
Universidad Michoacana\\
58060, Morelia, Mich., M\'exico.\\
\hbox{\small\tt rcampos@zeus.umich.mx, claudio@fismat.umich.mx}\\
}
\date{}
\maketitle
{\vskip.3cm
\noindent MSC: 65D25, 41A05, 42A15
}\\
\vspace*{1.5truecm}

\noindent{\large\bf Abstract}
\vspace*{.5cm}
\\
A procedure to obtain differentiation matrices is extended straightforwardly to yield new differentiation matrices useful to obtain 
derivatives of complex rational functions. Such matrices can be used to obtain numerical solutions of some singular differential problems defined in the complex domain. The potential use of these matrices is illustrated with the case of elliptic functions.
\vskip1.5cm

\section{Introduction}
In a sequel of papers (\cite{Cam97}-\cite{Cam04} and references therein),
a Galerkin-type method has been used to solve boundary value problems
and to find limit-cycles of nonautonomous dynamical systems. The method
is based on the discretization of the differential problem by using differentiation
matrices that are projections of the derivative in spaces of algebraic polynomials 
or trigonometric polynomials. This kind of matrices arise naturally in the 
context of interpolation of functions and yield exact values for the derivative of 
polynomial functions at certain points selected as nodes. 
The class of functions that defines the domain of the differential operator  
determines the kind of differentiation matrix to be used. Thus, to get a discrete form of 
an operator acting on real functions that drop off rapidly to zero at large distances, we can 
use the skew-symmetric differentiation matrix 
\begin{equation}\label{eq1}
D_{jk}=\begin{cases}0,&i=j,\\\noalign{\vskip .5truecm}
\displaystyle {{(-1)^{j+k}}\over{x_j-x_k}}, &i\not=j,\\ 
\end{cases}
\end{equation}
constructed with the $N$ zeros $x_j$ of the Hermite polynomial $H_N(x)$ as lattice 
points\footnote{Equation (\ref{eq1}) gives an asymptotic expression for $D_x$.}. 
On the other hand, to get a discrete form of an operator acting on periodic functions \cite{Cam01}, the 
differentiation matrix
\begin{equation}\label{eq2}
D_{jk}=\begin{cases}0,&j=k,\\\noalign{\vskip .5truecm}
\displaystyle {(-1)^{j+k}\over{2\sin{{(x_j-x_k)}\over 2} }}, &j\not=k,\\
\end{cases}
\end{equation}
should be used. In this case the lattice points can be chosen as the $N$ equidistant 
points
\begin{equation}
x_j=-\pi+2\pi j/ N, \quad j=1,2,\cdots,N.
\label{eq3}
\end{equation}
Discrete forms of multidimensional differential operators can be obtained by using 
direct products of suitable one-dimensional differentiation matrices (see \cite{Cam00} 
for instance).\\
Other approaches to obtain differentiation matrices can be found in \cite{Wel97}.\\
The purpose of this paper is to use the complex Hermite interpolation formula to find 
new differentiation matrices on the complex domain that produce an exact formula for 
the derivatives of complex rational functions and that can be used to approximate 
the solution of a singular differential equation in the complex domain. As examples, 
the derivatives of meromorphic two-periodic functions such as Jacobi elliptic functions 
and Weierstrass elliptic functions will be used to show the potential of this technique.\\


\section{Interpolatory differentiation matrices}\label{MatricesdeDif}
In order to illustrate the procedure to generate differentiation matrices, an already known 
matrix \cite{Cal81}, which gives the derivative of algebraic polynomials of the complex 
variable $z$ is obtained by using the Hermite interpolation formula.\\
\subsection{Algebraic polynomials}
Let $f(z)$ be an analytic function on a domain $G$ containing a closed rectifiable Jordan
curve $\gamma$ and let $z_k$ be $N$ different points of $I(\gamma)$ defining the polynomial 
\[
\omega(z)=\prod^N_{k=1}(z-z_k).
\]
Then, the unique polynomial $p(z)$ interpolating $f(z)$ associated to the set of 
points $z_k$ is given by (\cite{Mar85},  \cite{Wal69})
\begin{equation}
p(z)=\frac{1}{2\pi i}\int_\gamma \frac{f(\zeta)}{\omega(\zeta)} \frac{\omega(\zeta)-\omega(z)}{\zeta-z} d\zeta .
\label{eq4}
\end{equation}
The residual function $R(z)=f(z)-p(z)$, $z\in I(\gamma)$, vanishes at $z_k$ yielding that $f(z_k)=p(z_k)$. 
To see this, write $R(z)$ as
\begin{equation}
R(z)=\frac{\omega(z)}{2\pi i}\int_\gamma \frac{f(\zeta)}{(\zeta-z)\omega(\zeta)} d\zeta .
\label{eq4p1}
\end{equation}
Since the integral of the right-hand side of this equation represents an analytic function on $I(\gamma)$
we have that $R(z_k)=0$. 
\\
To obtain the form of $p(z)$, the integral in Eq. (\ref{eq4}) can be calculated by the residue theorem.
Since $\omega(\zeta)-\omega(z)$ is divisible by $\zeta-z$, the integrand has simple poles only at 
$z_1,\ldots,z_N$. Thus, the residue theorem yields the well-known Lagrange interpolation  formula
\begin{equation}
p(z)=\sum^N_{j=1} f(z_j) {{\prod_{k\neq j}^N (z-z_k)}\over{\prod_{k\neq j}^N (z_j-z_k)}}, 
\qquad z\in G.
\label{eq5}
\end{equation}
The differentiation matrix $D$ associated to this formula can be obtained if the derivative of 
(\ref{eq5}) is evaluated at $z_i$ and written as  
\begin{equation}
\frac{dp(z_i)}{dz} = \sum^N_{j = 1} D_{ij}\,f(z_j).
\label{eq6}
\end{equation}
Since
\[
\frac{d}{dz}\left[\prod_{k\neq j}^N (z-z_k)/\prod_{k\neq j}^N (z_j-z_k)\right] _{z=z_i} =
\begin{cases} \displaystyle{{\prod_{k \neq i,j}^N (z_i-z_k) }\over{\prod_{k \neq j}^N (z_j-z_k)}},&
 i\neq j, \\\noalign{\vskip .5truecm}
\displaystyle{\sum_{k \neq i}^N \frac{1}{(z_i-z_k)}},&i= j,\\
\end{cases}
\]
we get that the differentiation matrix for algebraic polynomials is given by
\begin{equation}\label{eq7}
D_{ij}=\begin{cases}\displaystyle \sum_{k \neq i}^N \frac{1}{(z_i-z_k)},&i=j,\\\noalign{\vskip .5truecm}
\displaystyle \frac{\omega '(z_i)}{(z_i-z_j)\,\,\omega '(z_j)}, &i\not=j.\\
\end{cases}
\end{equation}
If $f(z)$ is a polynomial of degree at most $N-1$, $f(z)$ is given identically by $p(z)$ and therefore 
formula (\ref{eq6}) gives the exact derivative of $f(z)$. Since $f'(z)$ is another polynomial of this class,
the derivatives of higher order can be obtained by applying successively the matrix $D$ to the vector 
of values $f(z_j)$, i.e.,
\begin{equation}
f^{(n)}=D^n f, \qquad n=0,1,2,\ldots.
\label{padnp}
\end{equation}
Here, $f^{(n)}$ is the vector whose entries are $d^nf(z_i)/dz^n$, $D^n$ is the $n$th power of $D$ and 
$f$ is just the vector whose entries are $f(z_i)$. The functional form of $f^{(n)}(z)$ can be obtained through 
an interpolation of the values yielded by (\ref{padnp}). Since any set of $N$ different complex
numbers belonging to $G$ yields the same polynomial $f(z)$, this result is independent of the points 
$z_k$. On the other hand, if $f(z)$ is not a polynomial 
of degree at most $N-1$ a residual vector should be added to the right-hand side of (\ref{eq6}) to
get $f'(z_i)$. However, such a formula yields an good approximation to $f'(z_i)$ 
if the absolute value of the $M$-th term of the Taylor series of $f(z)$ goes to zero rapidly as 
$M$ is increased.\\
\subsection{Trigonometric polynomials}
The preceding arguments can be modified to consider the interpolation of periodic 
functions in terms of trigonometric polynomials. Let $f(z)$ be a one-periodic analytic function with 
period $2\pi$ and let $G$ be a domain of
the open strip $-\pi < \Re z <\pi$, $-\infty < \Im z <\infty$, containing a closed rectifiable 
Jordan curve $\gamma$. \\
Since any trigonometric polynomial  $\tau(z)=a_0 + \sum^m_{k=1} (a_k \cos kz +b_k \sin kz)$ of 
degree at most $m$ can be written in the form $\tilde{\tau}(s)=s^{-m} q(s)$ under the the mapping 
$s=\varphi(z)= e^{iz}$ where $q(s)$ is a polynomial of degree at most $2m$ in $s$, we need to take an 
odd number $N=2m+1$ of different points $s_k\in \varphi(I(\gamma))$, i.e, $2m+1$ different points 
$z_k\in I(\gamma)$, to yield an exact interpolation formula in the case in which $f(z)$ is a 
trigonometric polynomial of degree at most $m$. This fact can be shown as follows.\\
Let us take $N=2m+1$. The set of points $s_k$, $k=1,2,\ldots ,N$ define the polynomial 
$\tilde{\omega}(s)=\prod^{N}_{k=1}(s-s_k)$. The interpolant function 
$\tilde{p}(s)$ to $\tilde{f}(s)=f(\varphi^{-1}(s))$ corresponding to the set of $N$ points $s_k$ is given by
\begin{equation}
\tilde{p}(s)=\frac{s^{-m}}{2\pi i}\int_{\tilde{\gamma}} \frac{\tilde{f}(\zeta)}{\tilde{\omega}(\zeta)} 
\frac{ s^{m} \tilde{\omega}(\zeta)- \zeta^{m}\tilde{\omega}(s)}{\zeta-s} d\zeta,
\label{eq8}
\end{equation}
where $\tilde{\gamma}=\varphi(\gamma)$. Since 
$[s^{m} \tilde{\omega}(\zeta)- \zeta^{m}\tilde{\omega}(s)]/(\zeta-s)$
is a polynomial in $s$ of degree $N-1=2m$, $\tilde{p}(s)$ has the required form $s^{-m} q(s)$,
where $q(s)$ is a polynomial of degree at most $2m$, to represent a trigonometric polynomial $\tau(z)$.\\
To show that $\tilde{f}(s_k)=\tilde{p}(s_k)$, let us consider the residual function 
$\tilde{R}(s)=\tilde{f}(s)-\tilde{p}(s)$ which is now
\[
\tilde{R}(s)=\frac{1}{2\pi i}\frac{\tilde{\omega}(s)}{s^m}\int_{\tilde{\gamma}} \frac{\tilde{f}(\zeta)\zeta^m}
{(\zeta-s)\tilde{\omega}(\zeta)} d\zeta. 
\]
By definition, $\tilde{G}=\varphi(G)$ does not contains points $s_k (\hbox{mod}\,\, 2\pi)$ other than 
$s_k$ therefore, the integral of the right-hand side of this equation represents an analytic function in 
$I(\tilde{\gamma})$ and we have that $\tilde{R}(s_k)=0$. \\
Since $s^{m} \tilde{\omega}(\zeta)- \zeta^{m}\tilde{\omega}(s)$ is divisible by  $\zeta-s$, the poles of 
the integrand are simple and located at $s_k$. The residue theorem yields now
\begin{equation}
\tilde{p}(s)=\sum^{N}_{j=1} \tilde{f}(s_j) \left(\frac{s_j}{s}\right)^m
{{\prod^{N}_{k\neq j} (s-s_k)}\over{\prod^{N}_{k\neq j} (s_j-s_k)}}, \qquad s\in \tilde{G}
\label{eq9}
\end{equation}
and the trigonometric polynomial of degree $m=(N-1)/2$ that interpolates to $f(z)$ is 
\begin{equation}
p(z)=\sum^{N}_{j=1}\displaystyle{ f(z_j) e^{i (N-1)(z_j-z)/2}
{{\prod^{N}_{k\neq j} (e^{i z}-e^{i z_k})}\over{\prod^{N}_{k\neq j} (e^{i z_j}-e^{i z_k})}} }, \qquad z\in G.
\label{eq10}
\end{equation}
Since 
\[
\displaystyle{e^{i (N-1)(z_j-z)/2}{{\prod^{N}_{k\neq j} (e^{i z}-e^{i z_k})}\over{\prod^{N}_{k\neq j} (e^{i z_j}-e^{i z_k})}}}=
\displaystyle{{\prod^{N}_{k\neq j} \sin\left( \frac{z-z_k}{2}\right) }\over{ \prod^{N}_{k\neq j} 
\sin\left(\frac{z_j-z_k}{2}\right)}},
\]
we obtain from (\ref{eq10}) the Gauss interpolation formula
\begin{equation}
p(z)=\sum^{N}_{j=1}f(z_j) 
{{\prod^{N}_{k\neq j} \sin\left(\frac{z-z_k}{2}\right)}\over{\prod^{N}_{k\neq j} \sin\left(\frac{z_j-z_k}{2}\right)}}, 
\quad N=2m+1, \quad z\in G.
\label{eq11}
\end{equation}\\
By writing the derivative of this formula in the form given by (\ref{eq6}) we can obtain the differentiation
matrix for trigonometric polynomials. Thus, the matrix $D$ whose entries are given by
\begin{equation}\label{dmtp}
D_{ij}=\begin{cases}\displaystyle \frac{1}{2}\sum_{k\neq j}^N\cot\left(\frac{z_j-z_k}{2}\right), &i=j,\\
\noalign{\vskip .5truecm}
\displaystyle \frac{1}{2}\csc\left(\frac{z_i-z_j}{2}\right)
\frac{\prod_{k\neq i}^N\sin(\displaystyle{\frac{z_i-z_k}{2})}}
{\prod_{k\neq j}^N\sin(\displaystyle{\frac{z_j-z_k}{2})}}, &i\neq j,\\
\end{cases}
\end{equation}
in terms of $N=2m+1$ different points $z_k\in G$, is a projection of $d/dz$ in the subspace of
trigonometric polynomials of degree at most $(N-1)/2$. Therefore, if $f(z)$ belongs to this space, 
$p(z)\equiv f(z)$, and $f^{(n)}(z)$ satisfies again an equation like (\ref{padnp}) but in this case
$D^n$ is the $n$th power of (\ref{dmtp}). The form of $f^{(n)}(z)$ can be obtained from the set of 
values $f^{(n)}(z_i)$ through an interpolation. For a general one-periodic analytic function with 
period $2\pi$ a residual vector should be added to (\ref{eq11}).\\ 
It is worth to notice that in the case in which the set of points $z_k$ are real numbers, the matrix (\ref{dmtp}) 
becomes the matrix used previously in \cite{Cam04}, \cite{Cam01}. However, if the points lay on a straight 
line which is parallel to the imaginary axis, $D$ takes the form
\[
D_{ij}=\begin{cases}-\displaystyle{\frac{i}{2}}\sum_{k\neq i}^N\coth\left(\frac{y_i-y_k}{2}\right), &i=j,\\
\noalign{\vskip .5truecm}
-\displaystyle{\frac{i}{2}}{\rm csch}\left(\frac{y_i-y_j}{2}\right)
\frac{\prod_{k\neq j}^N\displaystyle{\sinh(\frac{y_i-y_k}{2})}}{\prod_{k\neq i}^N
\displaystyle{\sinh(\frac{y_j-y_k}{2})}},
&i\neq j,\\ \end{cases}
\]
where $y_k= \Im z_k$, at the same time that the polynomial to differentiate becomes a linear combination 
of real hyperbolic functions.\\
\subsection{Rational functions}
Equation (\ref{eq8}) suggests the form of a interpolant rational function in the case
in which $f(z)$ is a meromorphic function.\\
Let $G$ be a domain that contains a closed rectifiable Jordan curve $\gamma$ and
let $f(z)$ be a meromorphic function with only one pole at $z=\alpha\not\in G$ of order $m$. 
Let us choose $N$ different points $z_k$ of $I(\gamma)$ and construct the polynomial 
$\omega(z)=\prod^N_{k=1}(z-z_k)$.  Thus, the rational function interpolating to $f(z)$
corresponding to the set of points $z_k$ is given by
\begin{equation}
p(z)=\frac{(z-\alpha)^{-m}}{2\pi i}\int_{\gamma} \frac{f(\zeta)}{\omega(\zeta)}\, 
\frac{ (z-\alpha)^{m} \omega(\zeta)- (\zeta-\alpha)^{m}\omega(z)}{\zeta-z} d\zeta.
\label{rfi}
\end{equation}
This can be shown by the same arguments used in the previous case. Again, 
$(z-\alpha)^{m} \omega(\zeta)- (\zeta-\alpha)^{m}\omega(z)$ 
is divisible by $\zeta-z$. Let $K^m_N(z,\zeta)$ be such a quotient, i.e.,
\[
K^m_N(z,\zeta)=\frac{ (z-\alpha)^{m} \omega(\zeta)- (\zeta-\alpha)^{m}\omega(z)}{\zeta-z}.
\]
Since $K^m_N(z,\zeta)$ is a polynomial of degree $N-1$ in $z$, $p(z)$ is a rational function of
form $q(z)/(z-\alpha)^{m}$ that interpolates to $f(z)$ at $z_k$, where $q(z)$ is a polynomial of 
degree at most $N-1$. The residual function
\[
R(z)=\frac{1}{2\pi i}\frac{\omega(z)}{(z-\alpha)^m}\int_{\gamma} \frac{f(\zeta)
(\zeta-\alpha)^m}{(\zeta-z)\omega(\zeta)} d\zeta. 
\]
vanish at $z_k$ because the integral is an analytic function in $I(\gamma)$ and
$\alpha\not\in G$. Therefore, $p(z_k)=f(z_k)$. The residue theorem yields now
\begin{equation}
p(z)=\sum^{N}_{j=1} f(z_j) \left(\frac{z_j-\alpha}{z-\alpha}\right)^m
{{\prod^{N}_{k\neq j} (z-z_k)}\over{\prod^{N}_{k\neq j} (z_j-z_k)}}, \qquad z\in G.
\label{rfrt}
\end{equation}
The derivative of this equation at $z_i$ can be written in the form (\ref{eq6}) where we have now
\begin{equation}
D_{ij}=\begin{cases}\displaystyle \sum_{k \neq i}^N \frac{1}{(z_i-z_k)}-\frac{m}{z_i-\alpha},
&i=j,\\\noalign{\vskip .5truecm}
\displaystyle \frac{(z_j-\alpha)^m/(z_i-\alpha)^m}{z_i-z_j}\,
\frac{\omega '(z_i)}{\omega '(z_j)}, &i\not=j.\\\end{cases}
\label{mdrf1}
\end{equation}
Obviously, if $f(z)$ is a rational function of the form $q(z)/(z-\alpha)^{m}$ where $q(z)$ is a 
polynomial of degree at most $N-1$, $f(z)\equiv p(z)$ and formula (\ref{eq6}) becomes
\begin{equation}
\frac{df(z_i)}{dz} = \sum^N_{j = 1} D_{ij}\,f(z_j). 
\label{dprf}
\end{equation}
However, the powers of $D$ do not give the derivatives of higher order as in the previous 
cases since $f'(z)$ does not has the required form: it has a pole of order $m+1$ at $z=\alpha$. 
Despite this, it is possible to obtain the value of $f^{(n)}(z)$ at $z_k$ by using the matrix 
\begin{equation}
D_n(m)=D_{m+n-1}D_{m+n-2}\cdots D_m.
\label{mdnm}
\end{equation}
Each matrix $D_k$ is defined by (\ref{mdrf1}) where the parameter $m$, defining the order of 
the pole, is substituted by each value of the index $k$. Thus, it should be clear that if 
$f(z)$ has the form given above, (\ref{padnp}) becomes 
\begin{equation}
f^{(n)}=D_n(m) f, \qquad n=0,1,2,\ldots.
\label{fnmrf}
\end{equation}
This equation can be generalized to the case in which $f(z)$ is a rational function of form
\begin{equation}
f(z)=\frac{\sum_{k=0}^{M}a_k z^k}{(z-\alpha_1)^{m_1}(z-\alpha_2)^{m_2}\cdots (z-\alpha_r)^{m_r} }
\label{frfr}
\end{equation}
where $\alpha_l\not\in G$, $l=1,2,\cdots,r$. A straightforward calculation gives the generalized 
form of (\ref{fnmrf})
\begin{equation}
f^{(n)}=D_n(m_1,m_2,\ldots,m_r) f, \qquad n=0,1,2,\ldots,
\label{gfnmrf}
\end{equation}
where now $D_n(m_1,m_2,\ldots,m_r)$ stands for the ordered matrix product 
\begin{equation}
D_n(m_1,m_2,\ldots,m_r)=\prod_{k=n-1}^0 D_{m_1+k,m_2+k,\ldots,m_r+k}
\label{dmrr}
\end{equation}
and $D_{\mu_1,\mu_2,\ldots,\mu_r}$ is the matrix whose entries are given by
\begin{equation}
(D_{\mu_1,\mu_2,\ldots,\mu_r})_{ij}=\begin{cases}\displaystyle \sum_{k \neq i}^N
 \frac{1}{(z_i-z_k)}-\sum_{k=1}^r\frac{\mu_k}{z_i-\alpha_k},
&i=j,\cr\noalign{\vskip .5truecm}
\displaystyle\frac{1}{z_i-z_j}\,\,
\frac{\omega '(z_i)}{\omega '(z_j)}
\prod_{k=1}^r\left(\frac{ z_j-\alpha_k }{ z_i-\alpha_k }\right)^{\mu_k}
, &i\not=j.\\\end{cases}
\label{mdrfg}
\end{equation}
It should be noticed that (\ref{gfnmrf}) is an exact formula whenever $N\ge M+n r-1$, where $M$
is the degree of the polynomial in the numerator of (\ref{frfr}). The reason is that after each 
differentiation the numerator of the derivatives of $f(z)$ is a polynomial whose degree grows
by $r$. If the function $f(z)$ to differentiate has poles at $\alpha_1,\ldots, \alpha_r$ of orders 
$n_1,\ldots,n_r$ instead $m_1,\ldots,m_r$,  with $n_1<m_1,\ldots,n_r<m_r$,  formula (\ref{gfnmrf}) 
is still exact whenever $N\ge M+n r-1+\sum_{k=1}^r (m_k-n_k)$. The reason for this is that in such a case
$f(z)$ can be writen in the form (\ref{frfr}) where the numerator is now a polynomial of degree 
$M+\sum_{k=1}^r (m_k-n_k)$. Obviously, in this case is much better to use the differentiation matrix  
$D_n(n_1,n_2,\ldots,n_r)$ instead $D_n(m_1,m_2,\ldots,m_r)$ for numerical purposes.
\vskip1cm

\section{Applications}\label{Aplicaciones}
As stated before, the numerical solution of differential problems can be accomplished by the use of
differentiation matrices, and in the case of a differential problem in the complex domain, the differentiation
matrices introduced in this paper may be useful. To illustrate the potential of their use, we choose two
meromorphic cases which are important in applications: Jacobi elliptic functions and Weierstrass elliptic 
functions. In both cases it is possible to establish the numerical convergence of the results since the 
derivatives of these functions are known. We also obtain approximate solutions of a differential equation
with a regular singularity at $z=0$.
Before beginning these examples is convenient to alert the
reader to the fact that the numerical implementation of the matrices for rational functions given above
may need a high-precision code: in most cases, the usual 16-digit precision is not enough to obtain
accurate results.\\ 

\subsection{A rational function}
Let us consider the function
\begin{equation}\label{funrac}
f(z)=\frac{z^7+z+1}{z^{10}}.
\end{equation}
According to the results of the last section, to obtain the exact value of the $n$-th derivative of (\ref{funrac}),
[cf. Eq. (\ref{fnmrf})], it is necessary to choose $N > 7$ different points $z_k\ne 0$ in the complex plane to 
build the matrix (\ref{mdnm}) where obviously, $\alpha=0$ and $m=10$.\\
As an example, let us take the third derivative of (\ref{funrac}). Table \ref{tabla} shows the max-norm of the vector whose 
entries are $\vert[f^{(3)}(z_j)-\sum_{k=1}^N (D_3(10))_{jk}f(z_k)]/f^{(3)}(z_j)\vert$ (the relative error) in terms of $N$. The nodes 
were chosen to be as evenly spaced on the ray $z=(1+i)t$, $1/2< t \le 1$, i.e., $z_k=(1+i)(1+k/N)/2 $, $k=1,\ldots,N$, and the computations were made with the standard 16-digit precision. As it can be seen from Table \ref{tabla}, the error vanishes for values of $N$ greater than $7$, as expected.\\
\begin{table}
\begin{center}
\begin{minipage}{12cm}
\caption
{\small The norm of the relative error $E_N=\vert (f'''-p''')/f'''\vert_\infty$ for (\ref{funrac}) in terms of the number of nodes $N$.}
\end{minipage}
\end{center}
\vskip.2cm
\label{tabla}
\hbox to \textwidth{\hfill
\begin{tabular}{|c|c|} \hline
\raise 15pt\hbox{\null}\hskip 5pt $N$ \hskip 5pt & $\hskip 1cm E_N\hskip 1cm $\hskip 1cm \\
\hline{\raise 15pt\hbox{\kern-.4em}}
4 & 0.657\\
5 & 0.136 \\
6 & 0.155 $\times 10^{-1}$ \\
7 & 0.742 $\times 10^{-3}$ \\
8 & 0. $\times 10^{-16}$ \\
9 & 0. $\times 10^{-16}$ \\
10 & 0. $\times 10^{-16}$ \\
11 & 0. $\times 10^{-16}$ \\\hline
\end{tabular}
\hfill}
\vskip1cm
\end{table}
\subsection{Elliptic functions}
As is known, an elliptic function is a doubly periodic function which is analytic except at poles 
\cite{Bow61}-\cite{Abr72} and one of the two simplest cases of elliptic functions corresponds to Jacobi's
functions; the other corresponds to Weierstrass' $\mathcal{P}\hbox{-function}$.\\
The differentiation matrices given above do not apply to functions like these, however it is possible to build
a matrix which is expected to provide approximate values of $f^{(n)}(z_k)$ along straight lines inside the
fundamental paralelograms and near to the poles. Since an elliptic function $f(z)$ becomes
a one-periodic function if $z$ is constrained to move along a straight line defined by one of the periods, such 
a matrix can be constructed by using trigonometric polynomials 
divided by an algebraic polynomial with zeros of suitable orders taken at the poles of $f(z)$. This procedure
is equivalent to taking apart the matrix (\ref{mdrfg}) and incorporating only the singular terms in (\ref{dmtp}) to 
yield the new matrices
\begin{equation}
\tilde{D}_n(m_1,m_2,\ldots,m_r)=\prod_{k=n-1}^0 \tilde{D}_{m_1+k,m_2+k,\ldots,m_r+k},
\label{dmrrjac}
\end{equation}
where 
\begin{equation}
(\tilde{D}_{\mu_1,\mu_2,\ldots,\mu_r})_{ij}=\begin{cases}\displaystyle \sum_{k \neq i}^N
\cot\left(\frac{z_j-z_k}{2}\right) -\sum_{k=1}^r\frac{\mu_k}{z_i-\alpha_k},
&i=j,\cr\noalign{\vskip .5truecm}
\displaystyle \frac{1}{2}\csc\left(\frac{z_i-z_j}{2}\right)
\prod_{k\neq i}^N\frac{\sin(z_i-z_k)/2}
{\sin(z_j-z_k)/2 }
\prod_{k=1}^r\left(\frac{ z_j-\alpha_k }{ z_i-\alpha_k }\right)^{\mu_k}, &i\not=j.\\\end{cases}
\label{mdjac}
\end{equation}
To test the numerical performance of this matrix we choose two numerical examples. The
first one corresponds to Jacobi's function $f(z)=\hbox{sn}(z| \, \frac{1}{2}\,)$ which has two 
periods $4K$ and $2iK'$ and two simple poles at $iK'$ and $2K+iK'$, where
\[ 
K=\int_0^{\pi/2} \frac{d\theta}{\sqrt{1-(\sin^2\theta) /2}},\qquad K'=K\approx 1.854.
\]
Therefore, to build the differentiation matrix we take $2\pi$-periodic trigonometric polynomials, $r=2$, 
$\mu_1=\mu_2=1$, $\alpha_1=iK'$, 
and $\alpha_2=2K+iK'$ in (\ref{mdjac}), and to measure the approximation of $f'(z_j)$ by 
$p\,'(z_j)=\sum_{k=1}^N (\tilde{D}_1(1,1))_{jk}f(z_k)$ we use the max norm. Here,
$f'(z)=\hbox{cn}(z| \, \frac{1}{2}\,)\, \hbox{dn}(z| \, \frac{1}{2}\,)$.
The results are displayed in Figure 1, where the max-norm of the error 
$f'(z_j)-p\,'(z_j)$ is plotted against the number of nodes showing numerical convergence. The number 
of digits of precision used in the calculations is 16 and the nodes are chosen to be as evenly spaced on the 
ray $z=(2+i)t$, $1/2< t \le 1$, i.e., $z_k=(2+i)(1+k/N)/2 $, $k=1,\ldots,N$. The norm of the error  is $1.6\times 10^{-4}$ 
for $N=10$ 
and  $1.0\times 10^{-8}$ for $N=20$.\\
Our second example is the Weierstrass $\mathcal{P}\hbox{-function}$, which has a double pole at 
the origin and two periods $\omega_1$ and $\omega_2$. The derivative of $\mathcal{P}(z)$ satisfies
\[  (\mathcal{P}'(z))^2=4(\mathcal{P}(z))^3-g_2\mathcal{P}(z)-g_3,  \]
where $g_2$ and $g_3$ are the elliptic invariants which are related to both $\omega_1$ and 
$\omega_2$ \cite{Abr72}.\\
To approximate the first derivative of $\mathcal{P}(z)$ at the nodes, we need to take $n=1$ and 
$r=1$ in (\ref{dmrrjac}) and $\mu_1=2$ and $\alpha_1=0$ in (\ref{mdjac}). The precision of the 
calculations, the nodes and the kind of trigonometric polynomials used to construct the differentiation matrix are the same as above. The numerical results are displayed in Figure 1 
and show again numerical convergence with a small number of nodes. For $N=10$ the norm of the
error is $1.5\times 10^{-5}$ and $1.0\times 10^{-8}$ for $N=20$.
\vskip1cm
\hbox to \textwidth{\hfill\scalebox{0.8}{\includegraphics{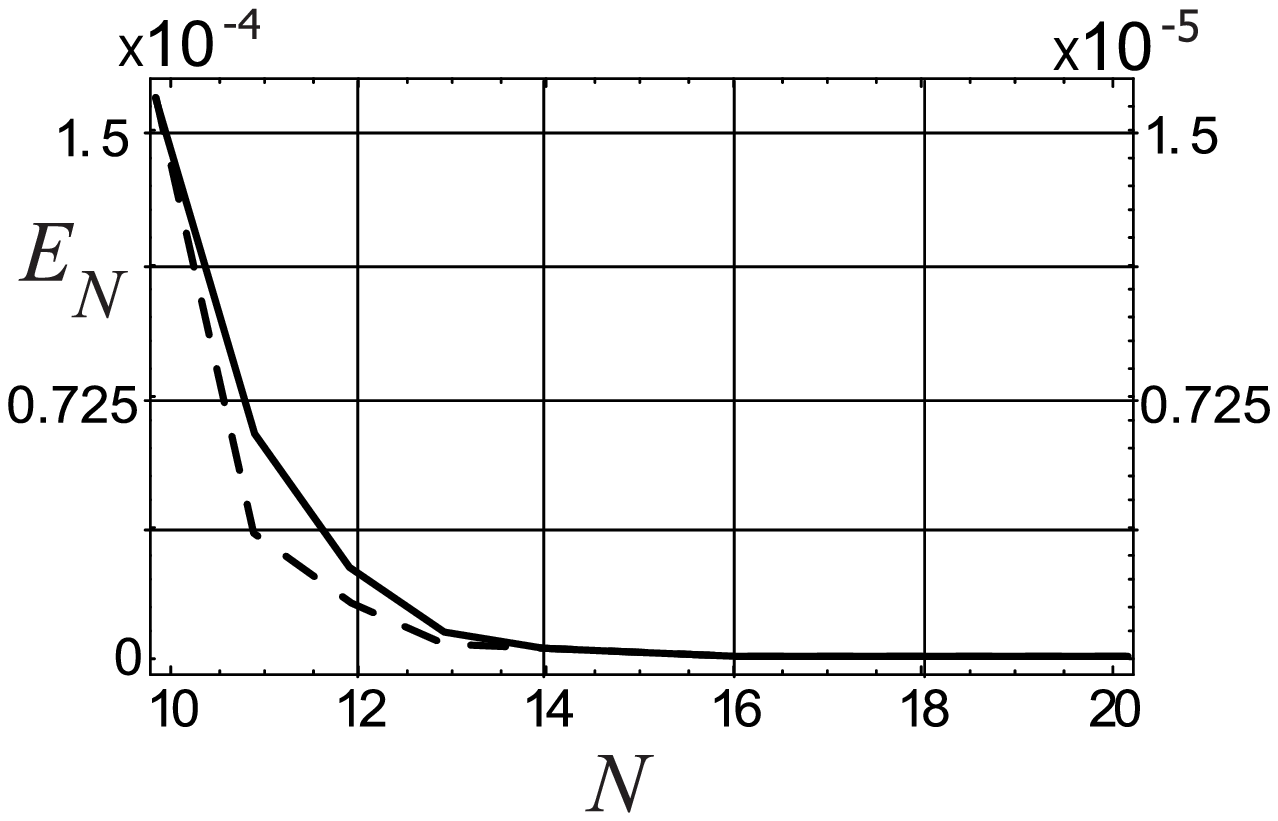}}\hfill}
\begin{center}
\begin{minipage}{14cm}
{\small 
Figure 1: The norm of the error $E_N=\vert f'-p'\vert_\infty$ against the number of nodes for the 
Jacobian elliptic function (solid line) and Weierstrass' $\mathcal{P}\hbox{-function}$ (broken
line).  The vertical axes are scaled accordingly to Jacobi's case (the left-hand side axis) or
Weierstrass' case (the right-hand side axis).
}
\end{minipage}
\end{center}
\subsection{Krummer's equation}
As a final example, we consider the Krummer differential equation, written as an eigenvalue problem
\begin{equation}\label{edokrum}
z\frac{d^2f(z)}{dz^2}+(b-z)\frac{df(z)}{dz}=a f(z),
\end{equation}
which has a regular singularity at $z=0$ and an irregular singularity at $\infty$. As is well known, the single-valued solution of this
equation is the confluent Hypergeometric function $M(a,b,z)=_1\hskip -4ptF_1(a,b,z)$. Since this function can be approximated by algebraic polynomials for $b\ne -n$ ($n$ a positive integer), we can obtain approximate solutions of this differential equation by using the differentiation matrix $D$ given by (\ref{dmtp}) to approximate the derivative of $M(a,b,z)$ and solving the 
$N$-dimensional eigenvalue problem
\begin{equation}\label{edokrumdis}
Lf_\lambda= \lambda f_\lambda,\quad L=Z D^2 +(b 1_N-Z) D
\end{equation}
where $Z$ is a diagonal matrix whose nonzero elements are the nodes $z_1,\ldots, z_N$, $b$ is a complex number ($b\ne -n$),
$1_N$ is the identity matrix of dimension $N$ and the eigenvalue $\lambda$ is the value of $a$ at which $M(a,b,z_k)$ is to be approximated by $(f_\lambda)_k$.  Let us denote by $M_a$ the vector whose $k$th component is $M(a,b,z_k)$. Since the $n$th coefficient of the power series of $M(a,b,z)$ is given by
\[
\frac{a(a+1)(a+2)\ldots (a+n-1)}{b(b+1)(b+2)\ldots (b+n-1)n!},
\]
the best approximation obtained for a given set of parameters $b$, $N$,  
$z_k$, is given by the eigenvector $f_\lambda$ ($\lambda=a$) corresponding to the eigenvalue with lowest absolute value $\lambda_m$. To construct the matrix $L$ in (\ref{edokrumdis}) we choose $N=21$ and $z_k=5(1+i) k /N$, $k=1,\ldots N$.  For $b$ we take the cases $b=5/2$ and $b=3+2i$.  In order to compare the approximate and exact results, we normalize both vectors 
$M_\lambda{_m}$ and $f_\lambda{_m}$ with the max-norm. The calculations were made with 16 digits of precision and the results are displayed in Fig. 2. The absolute errror $\Vert M_\lambda{_m}-f_\lambda{_m} \Vert_0$ is $0.0675659$ for $b=5/2$ and $0.0426948$ for $b=3+2i$.
\vskip1cm
\hbox to \textwidth{\hfill\scalebox{0.8}{\includegraphics{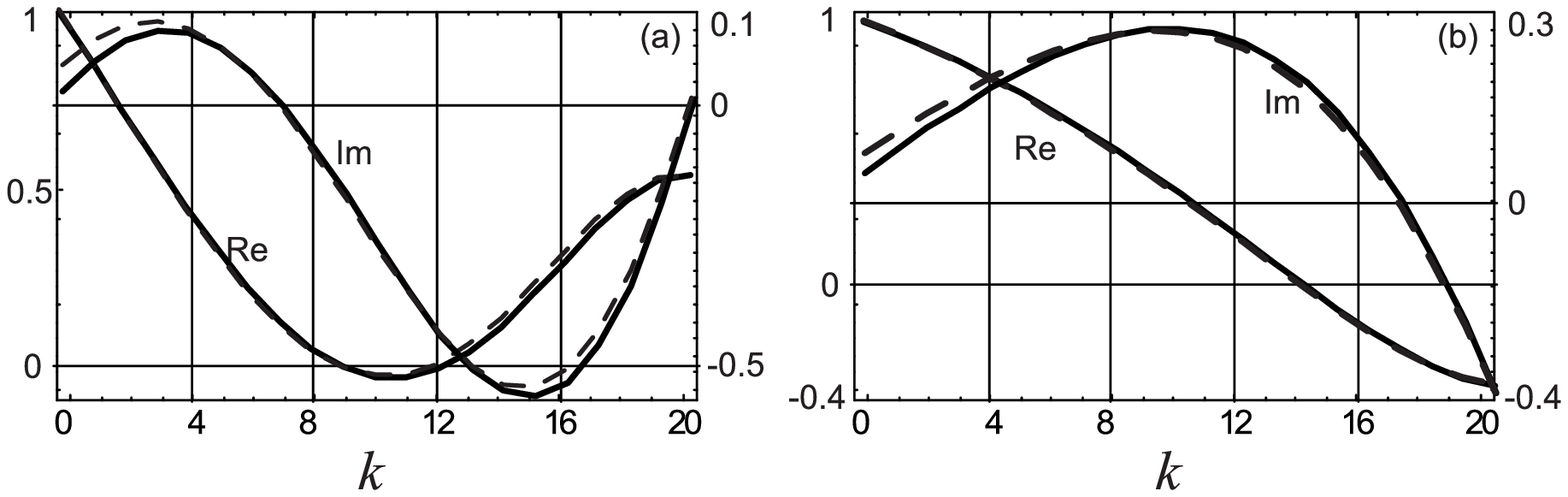}}\hfill}
\begin{center}
\begin{minipage}{14cm}
{\small 
Figure 2: Normalized real (Re) and imaginary (Im) parts of $f_\lambda{_m}$ plotted versus their index (solid lines). Case (a) corresponds to $b=5/2$ and $\lambda_m=-0.301513+1.00758 i$ and case (b) to $b=3+2i$ and  
$\lambda_m=-0.381925+0.527533 i$. They are compared with the exact values of the Krummer function (broken lines). The matrix $L$ is constructed with $21$ nodes $z_k=5(1+i) k /N$, $k=1,\ldots 21$. The scale on the left-hand side vertical axis corresponds to the real part and the one on the right-hand side to the imaginary part. 
}
\end{minipage}
\end{center}
\vskip1cm
\section{Concluding remark}
According to the results of section \ref{MatricesdeDif}, the process of interpolation in vector 
spaces of polynomials of dimension $N$ maps the derivative $d/dx$ into a $N\times N$ 
matrix $D$. The fact that a differential operator acting on a vector space of finite dimension 
can be written as a matrix is not a surprise of course, however it should be noted 
that the matrix $D$ yields the derivative of a function by taking the values of the function at 
$N$ arbitrary (but different) points including the point where the derivative is to be evaluated, 
i.e., it acts on a function as a nonlocal operator in spite of the local character of a differential
operator as the derivative.


\end{document}